\documentclass[11pt]{amsart}
\usepackage{amssymb}
\usepackage{amsmath}
\usepackage{mathtools}
\usepackage{enumitem}
\usepackage{mathrsfs}
\usepackage{hyperref}
\usepackage{comment}
\usepackage{tikz-cd}

\usepackage{stmaryrd}

\usepackage{fullpage}

\newtheorem{theorem}{Theorem}[section]
\newtheorem{lemma}[theorem]{Lemma}
\newtheorem{proposition}[theorem]{Proposition}
\newtheorem{corollary}[theorem]{Corollary}
\newtheorem{conjecture}[theorem]{Conjecture}

\newtheorem{alphtheorem}{Theorem}

\theoremstyle{definition}
\newtheorem*{ack}{Acknowledgments}

\newtheorem{remark}[theorem]{Remark}

\newtheorem{definition}[theorem]{Definition}
\newtheorem{question}[theorem]{Question}

\numberwithin{equation}{section} \numberwithin{figure}{section}

\DeclareMathOperator{\Pic}{Pic}

\DeclareMathOperator{\Sym}{Sym}

\DeclareMathOperator{\MW}{MW}
\DeclareMathOperator{\Num}{Num}

\DeclareMathOperator{\id}{id}

\newcommand*\ratmap{\mathbin{\tikz [baseline=0ex,-latex, dashed, ->] \draw [densely dashed] (0em,0.58ex) -- (1.3em,0.58ex);}}

\title{New examples of geometrically special varieties: K3 surfaces, Enriques surfaces, and algebraic groups}

\author{Finn Bartsch}
\address{Finn Bartsch \\
IMAPP Radboud University Nijmegen \\
PO Box 9010, 6500GL \\
Nijmegen, The Netherlands\\}
\email{f.bartsch@math.ru.nl}



\begin{document}

\begin{abstract}
We verify that elliptic K3 surfaces and algebraic groups have many rational points over function fields, i.e., they are geometrically special in the sense of Javanpeykar--Rousseau.
We also show that under additional assumptions, this geometric specialness persists under removal of closed subsets of codimension at least two.
\end{abstract}

\maketitle
\thispagestyle{empty}

\section{Introduction}

This article is concerned with the notion of a geometrically special variety as defined by Javanpeykar and Rousseau in \cite[Definition~1.7]{JRGeomSpec}.
Let $k$ be an algebraically closed field of characteristic zero.
A \emph{variety} is an integral separated finite type scheme over $k$.

\begin{definition}\label{def:maindefinition} 
Let $X$ be a variety.
We say that $X$ is \emph{geometrically special} if there is a Zariski-dense subset $S \subseteq X(k)$ such that for every $s \in S$ there is a smooth quasi-projective curve $C$, a point $c \in C$ and a collection of morphisms $(\phi_i \colon C \to X)_{i \in I}$ satisfying $\phi_i(c) = s$ such that $\bigcup_{i \in I} \Gamma_{\phi_i} \subseteq C \times X$ is Zariski-dense.
Here, $\Gamma_{\phi_i}$ denotes the graph of the morphism $\phi_i$.
\end{definition}

Slightly abusing the language, we will call a collection of morphisms $(\phi_i)_{i \in I}$ as in the above definition a \emph{covering set for $X$ through $s$}, even though what is really being covered is the product space $C \times X$.
A variety $X$ being geometrically special gives a precise meaning to the idea that ``$X$ has lots of rational points over some function field $k(C)$'' which is stronger than merely asking for the nonconstant $k(C)$-rational points to be Zariski-dense in $X$.

Our main results then can be summarized as follows:

\begin{alphtheorem}[Propositions \ref{prop:groupschemegeomspec} and \ref{prop:puncturedgroupschemegeomspec}] \label{mainthm_groups}
Let $X$ be a connected algebraic group.
Then $X$ is geometrically special.
If the base field $k$ is moreover assumed to have infinite transcendence degree over $\mathbb{Q}$ and $\Delta \subseteq X$ is any closed subset of codimension $\geq 2$, then $X \setminus \Delta$ is geometrically special as well. 
\end{alphtheorem}

\begin{alphtheorem}[Corollaries \ref{cor:ellk3geomspec} and \ref{cor:punctured_k3_geomspec}] \label{mainthm_k3}
Let $X$ be a K3 surface admitting an elliptic fibration.
Then $X$ is geometrically special.
If moreover $\Delta \subseteq X$ is a finite set of closed points in general position, then $X \setminus \Delta$ is geometrically special as well.
\end{alphtheorem}

\begin{alphtheorem}[Corollaries \ref{cor:enriquesgeomspec} and \ref{cor:punctured_enriques_geomspec}] \label{mainthm_enriques}
Let $X$ be an Enriques surface.
Then $X$ is geometrically special.
If moreover $\Delta \subseteq X$ is a finite set of closed points in general position, then $X \setminus \Delta$ is geometrically special as well.
\end{alphtheorem}

In Theorems~\ref{mainthm_k3}~and~\ref{mainthm_enriques}, we say that a property is satisfied for a subset in \emph{general position} if, for every integer $\delta \geq 0$, there is a dense open subset $U \subseteq X^{\delta}$ of the product space $X^{\delta} = X \times \cdots \times X$ such that the property is satisfied for every $\Delta \in U$.

\subsection{Motivation: Campana-special varieties}

In \cite[Conjecture~1.8]{JRGeomSpec}, inspired by Campana's conjectures \cite{CampanaFourier}, Javanpeykar and Rousseau conjecture the following:

\begin{conjecture}\label{conj:special}
Let $X$ be a smooth variety. Then $X$ is geometrically special if and only if $X$ is (Campana-)special.
\end{conjecture}

To understand this conjecture, we first need to introduce the class of (Campana-)special varieties, which were first defined by Campana in \cite{CampanaFourier} as a class ``opposite'' to the class of general type varieties.
To do so, let $X$ be a smooth proper variety equipped with a simple normal crossings divisor $D \subseteq X$.
We say that a line bundle $\mathcal{L}$ on $X$ is a \emph{Bogomolov sheaf for the pair $(X,D)$} if its Iitaka dimension $p := \kappa(\mathcal{L})$ is positive and there is a nonzero morphism $\mathcal{L} \to \Omega^p_X(\log D)$ into the sheaf of $p$-forms with at most logarithmic poles along $D$ (see \cite[Chapter~11]{IitakaBook} for the definition of the latter sheaf).
A smooth variety $X$ is then said to be \emph{Campana-special} if it has a proper compactification $\overline{X}$ whose boundary is a simple normal crossings divisor $D$ such that the pair $(\overline{X},D)$ has no Bogomolov sheaf.
This definition turns out to be independent of the chosen compactification \cite[Lemma~2.1]{Puncturing}.
It is then easy to check that a Campana-special variety does not admit a dominant morphism to a variety of log-general type -- indeed, if $f \colon X \to Y$ is such a morphism, saturating the line bundle $f^* \omega_Y$ gives rise to a Bogomolov sheaf on $X$.
This is the aforementioned sense in which Campana-special varieties are ``opposite'' to being of general type.

Conjecture~\ref{conj:special} was in part motivated by Lang's conjectured equivalence between the property of being of (log-)general type and various notions of hyperbolicity (see \cite{Lang86} or the survey article \cite[\S 12]{HyperbolicSurvey} for precise statements).
Indeed, it is easy to see that a geometrically special variety cannot dominate a variety which is algebraically hyperbolic in the sense of Demailly or more generally, geometrically hyperbolic (see \cite{HyperbolicSurvey} for the definitions of these notions).
In this sense, just like Campana-special varieties are ``opposite'' to being of general type, geometrically special varieties are ``opposite'' to geometrically hyperbolic varieties.
Hence, if we believe in Lang's conjectures that the class of algebraically (resp. geometrically) pseudo-hyperbolic varieties agrees with the class of varieties of log-general type, Conjecture~\ref{conj:special} appears to be rather reasonable. 
For more details on the motivation behind Conjecture~\ref{conj:special} and its relation to the analogous number theoretic conjecture dealing with potential density of integral points, we refer the reader to the introduction of \cite{JRGeomSpec}.

A smooth curve is Campana-special if and only if it is not of log-general type and it follows from the Theorem of de Franchis that this is the case if and only if the curve is geometrically special.
Consequently, Conjecture~\ref{conj:special} holds for curves.
In higher dimensions, examples of Campana-special varieties include all rationally connected varieties \cite[Theorem~3.22]{CampanaFourier} and all proper varieties of Kodaira dimension zero \cite[Theorem~5.1]{CampanaFourier}.
Another example is given by algebraic groups.
Since this appears to be missing from the literature, we give a proof of this fact in Lemma~\ref{lemma:alg_groups_special} below.
As K3 surfaces and Enriques surfaces are of Kodaira dimension zero, Theorems \ref{mainthm_groups}, \ref{mainthm_k3} and \ref{mainthm_enriques} above can be seen as partial confirmations of Conjecture~\ref{conj:special}.
Other cases in which Conjecture~\ref{conj:special} is known to be true are the case of closed subvarieties of abelian varieties and more generally, proper varieties of maximal Albanese dimension (where it follows from work of Kawamata, Ueno and Yamanoi, see \cite[Theorem~3.5 and Corollary~3.10]{JRGeomSpec} for details), and the case of symmetric powers of curves and those of certain surfaces (see \cite{Puncturing}).
Lastly, in \cite{PRTNumericallySpecial} it is shown that a smooth projective variety $X$ admitting a line bundle $\mathcal{L} \subseteq \Omega^1_X$ of \emph{numerical dimension} $1$ is neither Campana-special nor geometrically special.

\subsection{Motivation: ``Puncturing problems''}

If $X$ is a variety and $U \subseteq X$ is an open subvariety which is geometrically special, then $X$ is geometrically special as well.
Indeed, any covering set for $U$ can be postcomposed with the inclusion map $U \to X$ to obtain a covering set for $X$.
The converse is in general not true:
While $\mathbb{P}^1$ is easily seen to be geometrically special, its open subset $\mathbb{P}^1 \setminus \{0,1,\infty\}$ is not geometrically special.
To see this, we can for example appeal to Siegel's theorem over function fields which implies that for any smooth quasi-projective curve $C$, the morphism $(\mathbb{P}^1 \setminus \{0,1,\infty\}) \times C \to C$ has only finitely many sections whose projection to $\mathbb{P}^1 \setminus \{0,1,\infty\}$ is nonconstant.
Hence there are only finitely many nonconstant maps $C \to \mathbb{P}^1 \setminus \{0,1,\infty\}$.
However, one can still ask whether the converse holds if the codimension of $Z := X \setminus U$ in $X$ is at least two, i.e. one can ask: If $X$ is geometrically special and $Z \subseteq X$ is a closed subvariety of codimension $\geq 2$, is $X \setminus Z$ still geometrically special?
Indeed, in this case, one might expect that for a ``generic'' covering set $(\phi_i \colon C \to X)_{i \in I}$, ``most'' of the images $\phi_i(C)$ will be disjoint from the codimension $\geq 2$ subset $Z$ and hence that ``most'' covering sets for $X$ induce covering sets for $U$.
Nonetheless, phrased in this generality, the answer to this question is still negative.
A counterexample is given by $X = \Sym^g(C \times \mathbb{P}^1)$ with $C$ a smooth projective curve of genus $g \geq 2$ and $U \subseteq X$ the complement of the big diagonal; see \cite{Puncturing} for details.
However, this counterexample is a singular variety (albeit with canonical singularities), so that one might wonder whether this converse at least holds if $X$ is smooth (or has terminal singularities):

\begin{question}\label{q:puncturing}
Let $X$ be a smooth geometrically special variety and let $Z \subseteq X$ be a closed subvariety of codimension at least two. Is $X \setminus Z$ still geometrically special? What if $X$ is allowed to have terminal singularities?
\end{question}

This appears to be a difficult problem and the partial answers included in Theorems \ref{mainthm_groups}, \ref{mainthm_k3} and \ref{mainthm_enriques} are obtained using methods that are quite specific to the situation at hand, with rather little hope of generalizing them.

In this direction, let us also note that if $X$ is a smooth quasi-projective Campana-special variety and $Z \subseteq X$ is a closed subset of codimension at least two, the complement $X \setminus Z$ is again Campana-special \cite[Theorem~G]{Puncturing}.
Hence, Conjecture~\ref{conj:special} implies a positive answer to Question~\ref{q:puncturing} in the case when $X$ is smooth.

\begin{ack}
I thank Ariyan Javanpeykar for his constant support and many interesting and enlightening discussions.
I also thank Julian Demeio for a helpful discussion on Siegel's theorem over function fields, Remke Kloosterman for suggesting Lemma~\ref{lemma:numeq} and Manfred Lehn for suggesting the proof of Lemma~\ref{lemma:finitegeneratingset}.
\end{ack}

\section{Algebraic groups}

In this section, we show that for every algebraically closed field $k$ of characteristic zero, every connected algebraic group (i.e. finite type group scheme) over $k$ is geometrically special.
If in addition we also assume $k$ to have infinite transcendence degree over $\mathbb{Q}$, we moreover verify that the complement of a closed set of codimension at least two remains geometrically special.
Our main technical tool will be the following recently established criterion for testing density of sections in product spaces.

\begin{proposition}\label{prop:geomspeccriterion}
Let $Y$ be a variety and let $X$ be a quasi-projective variety. Let $(\phi_i \colon Y \to X)_{i \in I}$ be a family of morphisms. Suppose that there is a point $y_0 \in Y$ such that $\{ \phi_i(y_0)~|~i \in I\}$ is Zariski-dense in $X$. Then $S = \bigcup \Gamma_{\phi_i}$ is Zariski-dense in $Y \times X$.
\end{proposition}
\begin{proof}
See \cite[Theorem~D]{Puncturing}.
\end{proof}

We will also need the following well-known properties of algebraic groups.

\begin{lemma}[{\cite[Tag~0BF7]{Stacks}}] \label{lemma:groupschemequasiproj}
Every algebraic group over any field $K$ is quasi-projective over $K$.
\end{lemma}

\begin{lemma}[{\cite[Lemma~5.4]{JavAut}}] \label{lemma:infiniteorderpoint} 
Let $k$ be an algebraically closed field of characteristic zero and let $G$ be a nontrivial connected algebraic group over $k$. Then $G(k)$ contains an element of infinite order.
\end{lemma}

Given an algebraic group $G$ and two subsets $A, B \subseteq G(k)$, we denote by $A \cdot B$ the set consisting of all products $ab$ with $a \in A$ and $b \in B$.
We record the following elementary lemma about the behaviour of this operation with respect to taking closures.

\begin{lemma}\label{lemma:closure_alggroups}
Let $G$ be an algebraic group over $k$ and let $A, B \subseteq G(k)$ be any two subsets. Then $\overline{A \cdot B} = \overline{\overline{A} \cdot \overline{B}}$. 
\end{lemma}
\begin{proof}
The inclusion $\overline{A \cdot B} \subseteq \overline{\overline{A} \cdot \overline{B}}$ being obvious, we focus on the reverse inclusion. To show the reverse inclusion, it suffices to verify $\overline{A} \cdot \overline{B} \subseteq \overline{A \cdot B}$. We have
$$ \overline{A} \cdot \overline{B} = \bigcup_{a \in \overline{A}} a\overline{B} = \bigcup_{a \in \overline{A}} \overline{a B} $$
and thus, it suffices to verify $\overline{a B} \subseteq \overline{A \cdot B}$ for every $a \in \overline{A}$, for which it suffices to show $a B \subseteq \overline{A \cdot B}$ for every $a \in \overline{A}$. In other words, we need to show $\overline{A}\cdot B \subseteq \overline{A \cdot B}$. But now observe
$$ \overline{A} \cdot B = \bigcup_{b \in B} \overline{A} b = \bigcup_{b \in B} \overline{A b} $$
and the latter is clearly a subset of $\overline{A \cdot B}$. 
\end{proof}

If $G$ is an algebraic group and $g \in G(k)$ is an element such that the discrete subgroup $\langle g \rangle \subseteq G(k)$ generated by it is Zariski-dense, then we may prove geometric specialness of $G$ as follows:
Let $C \subseteq G$ be a curve containing $g$ and the neutral element $e$. Then the maps $C \to G$ obtained by post-composing the inclusion map with the $\alpha$-th power map $[\alpha] \colon G \to G$ for $\alpha \in \mathbb{Z}$ form a covering set for $G$ through the point $e$. Indeed, by construction, the images of $g$ under these countably many maps are a dense subset of $G$, whence the required density of the graphs of these maps in $C \times G$ follows by Proposition \ref{prop:geomspeccriterion}. If $G$ is an abelian variety, then an element $g \in G(k)$ as above always exists and the covering set constructed above is how the geometric specialness of abelian varieties is verified in \cite[Proposition~3.1]{JRGeomSpec} (although there, the final appeal to Proposition \ref{prop:geomspeccriterion} is replaced by an argument more specific to abelian varieties). However, for general algebraic groups $G$, there might not be an element $g$ for which $\langle g \rangle$ is Zariski-dense -- this happens already for $G = \mathbb{G}_a^2$. The following lemma provides a substitute.

\begin{lemma}\label{lemma:finitegeneratingset}
Let $k$ be an algebraically closed field of characteristic zero and let $G$ be a connected algebraic group over $k$.
Then there are finitely many elements $g_1,...,g_r \in G(k)$ such that the set
$$ \{ g_1^{\alpha_1}g_2^{\alpha_2}...g_r^{\alpha_r}~|~\alpha \in \mathbb{Z}^r \} $$
is Zariski-dense in $G$.
\end{lemma}
\begin{proof}
We first introduce some notation. If $g_1,...g_n \in G$ have already been chosen, we define for every $i=1,...,n$ the set $Z_i$ as the closure of the set $\{g_i^\alpha\}_{\alpha \in \mathbb{Z}}$ and
$$ Y_i := \{ g_1^{\alpha_1}g_2^{\alpha_2}...g^{\alpha_i}_i~|~\alpha \in \mathbb{Z}^i\}. $$
Observe that by Lemma \ref{lemma:closure_alggroups} we have $\overline{Y_i} = \overline{Z_1 \cdot Z_2 \cdot ... \cdot Z_i} = \overline{Y_{i-1}\cdot Z_i}$.

We now start the construction by choosing $g_1,...,g_n$ in such a way that the dimension of $\overline{Y_1},...,\overline{Y_n}$ is strictly increasing in each step until this is no longer possible.
At that point we stop and let $H$ be the connected component of the identity in $\overline{Y_n}$.
By Lemma \ref{lemma:infiniteorderpoint} we see that unless $G$ is trivial (in which case there is nothing to prove anyway), we have that $H$ is positive-dimensional.
We claim that $H$ is a closed normal subgroup scheme of $G$.

To see this, consider $\overline{H \cdot H}$.
Using Lemma \ref{lemma:closure_alggroups}, we see that $\overline{H \cdot H}$ is contained in $\overline{H \cdot Y_n} = \overline{H \cdot Z_1 \cdot ... \cdot Z_n}$.
Now note that $\overline{H \cdot Z_1}$ has to have the same dimension as $H$:
Otherwise, setting $g_{n+1}=g_1$ leads to $\overline{Y_{n+1}}$ having strictly bigger dimension than $\overline{Y_n}$, which is impossible by assumption.
Repeating this argument $n$ times, we see that $\overline{H \cdot Z_1 \cdot ... \cdot Z_n}$ has the same dimension as $H$.
So $\overline{H \cdot H}$ is a closed subset of $G$ which contains $H$ and has the same dimension as $H$.
Furthermore, there is a natural dominant map $H \times H \to \overline{H \cdot H}$ and since $H$ is connected, we see that $\overline{H \cdot H}$ is connected, too.
This implies $H \cdot H = H$.
The same arguments, applied to $H \cdot Z_n \cdot ... \cdot Z_1$, show that $H \cdot H^{-1} = H$, which implies $H^{-1} \subseteq H$ and thus that $H$ is a closed subgroup scheme of $G$.
Lastly, in the same way we also see that for every $g \in G(k)$, we have $H \cdot (g^{-1}Hg) = H$, which implies $g^{-1}Hg \subseteq H$.
Thus it follows that $H$ is normal, as desired.

Since we now know that $H \subseteq G$ is a closed normal subgroup scheme, we can form the quotient $G/H$.
It is a connected algebraic group of strictly lower dimension than $G$.
Doing an induction over $\dim(G)$, we may thus assume that the lemma has been proven for the group scheme $G/H$.
So we can find finitely many elements $h_1,...,h_s \in (G/H)(k)$ such that
$$ S := \{h_1^{\beta_1}...h_s^{\beta_s}~|~\beta\in\mathbb{Z}^s\} $$
is Zariski-dense in $G/H$.
Now choose any preimages $\tilde{h_1},...,\tilde{h_s} \in G(k)$. To complete the proof, it now suffices to show that the following set is dense in $G$:
$$ S' := \{g_1^{\alpha_1}...g_n^{\alpha_n}\tilde{h_1}{}^{\beta_1}...\tilde{h_s}{}^{\beta_s}~|~\alpha\in\mathbb{Z}^n, \beta\in\mathbb{Z}^s\} $$
From Lemma \ref{lemma:closure_alggroups} it follows that the closure of $S'$ contains all elements of the form $h\tilde{h_1}{}^{\beta_1}...\tilde{h_s}{}^{\beta_s}$ with $h \in H$.
Thus, $S'$ contains the preimage of $S$ under the quotient map $G \to G/H$.
Since $S \subseteq G/H$ is Zariski-dense and since quotient maps of algebraic groups are always flat, it follows that $S' \subseteq G$ is Zariski-dense.
This finishes the proof.
(It also follows that we had $H=G$ the entire time.) 
\end{proof}

\begin{remark}
We note that we can also use Lemma \ref{lemma:finitegeneratingset} to rederive the well-known fact that given an algebraic group $G$, its set of integral points is potentially dense. More precisely, we claim that if $G$ is defined over the ring of $S$-integers $\mathcal{O}_{K,S}$ of some number field $K$, there is a field extension $K \subseteq L$ and a set of places $S'$ of $L$ extending $S$ such that the $\mathcal{O}_{L,S'}$-points of $G$ are Zariski-dense. To see this, apply Lemma \ref{lemma:finitegeneratingset} to the base change $G_{\overline{\mathbb{Q}}}$ and obtain elements $g_1,...,g_r \in G(\overline{\mathbb{Q}})$. Now choose $L$ and $S'$ such that the finitely many points $g_1,...,g_r$ are defined over $\mathcal{O}_{L,S'}$. As the group law of $G$ is defined over $\mathcal{O}_{L,S'}$ as well, the products of the $g_i$ give rise to a Zariski-dense set of points defined over $\mathcal{O}_{L,S'}$.
\end{remark}

If we in addition assume our base field to be uncountable, we can furthermore arrange the choice of the elements $g_1,...,g_r$ in the previous lemma in such a way that the set of their products avoids any prespecified closed subset of codimension two.

\begin{lemma}\label{lemma:finitegeneratingset2}
Let $k$ be an uncountable algebraically closed field of characteristic zero and let $G$ be a connected algebraic group over $k$. Let $\Delta \subsetneq G$ be a proper closed subset not containing the origin. Then there are finitely many elements $g_1,...,g_r \in G(k)$ such that the set
$$ \{ g_1^{\alpha_1}g_2^{\alpha_2}...g_r^{\alpha_r}~|~\alpha \in \mathbb{Z}^r \} $$
is Zariski-dense in $G$ and disjoint from $\Delta$.
\end{lemma}
\begin{proof}
We reuse the notation $Y_i$ from the proof of Lemma \ref{lemma:finitegeneratingset}.
We claim that if we inductively choose very general elements $g_1,...,g_r$ until $\overline{Y_r} = G$, the condition that $Y_r$ is disjoint from $\Delta$ is automatically satisfied.
Here, \emph{very general} means that we have to avoid countably many proper Zariski-closed subsets of $G$ in each step (which is possible as we assumed $k$ to be uncountable). 

Indeed, suppose that for an integer $n \geq 1$, the elements $g_1,...,g_{n-1}$ were already chosen such that $Y_{n-1}$ is disjoint from $\Delta$ and such that $\overline{Y_{n-1}} \neq G$ (that this is possible for $n=1$, i.e. that $Y_0$ is disjoint from $\Delta$, is the assumption that $\Delta$ does not contain the origin).
Then we claim that a very general choice $g_n \in G(k)$ ensures that $\dim \overline{Y_n} > \dim \overline{Y_{n-1}}$ and $Y_n \cap \Delta = \emptyset$.
To see this, note that the condition $\dim \overline{Y_n} > \dim \overline{Y_{n-1}}$ simply requires the sets $(\overline{Y_{n-1}}g_n^\alpha)_{\alpha \in \mathbb{Z}}$ to be pairwise distinct.
This is given as soon as $g_n^{\alpha} \notin \overline{Y_{n-1}}$ for every $\alpha \in \mathbb{Z}$, i.e. $g_n \notin [\alpha]^{-1}\overline{Y_{n-1}}$ for every $\alpha \in \mathbb{Z}$, where $[\alpha]$ denotes the multiplication-by-$\alpha$ morphism $G \to G$.
As the $[\alpha]^{-1}\overline{Y_{n-1}}$ are countably many proper Zariski-closed subsets, this shows the desired claim.
For the condition that $Y_n \cap \Delta = \emptyset$, assume that we have an element $h \in Y_n \cap \Delta$.
Then, as $h \in Y_n$, we can write $h = h' g_n^{\alpha}$ for some $h' \in Y_{n-1}$ and some $\alpha \in \mathbb{Z}$.
This implies that $g_n \in [\alpha]^{-1}(h'^{-1} \Delta)$.
As there are only countably many choices for $h'$ and $\alpha$, respectively, we see that for a very general $g_n$, the intersection $Y_n \cap \Delta$ must be empty.
This finishes the proof.
\end{proof}

We can now prove that algebraic groups are geometrically special.

\begin{proposition}\label{prop:groupschemegeomspec}
Let $k$ be an algebraically closed field of characteristic zero.
Let $G$ be a connected algebraic group over $k$.
Then $G$ is geometrically special.
\end{proposition}
\begin{proof}
Let $g_1,...,g_r \in G(k)$ be elements as in Lemma \ref{lemma:finitegeneratingset}. Consider the variety $G^r = G \times ... \times G$. It is quasi-projective (Lemma \ref{lemma:groupschemequasiproj}). By Bertini's theorem, we find a smooth quasi-projective curve $C \subseteq G^r$ connecting $(e,e,e,...,e)$ (where $e \in G(k)$ is the neutral element) and $(g_1,...,g_r)$. Denote by $\iota$ the inclusion $C \to G^r$ and let $c_e$ and $c_g$ be the two points on $C$ just described. 

For every $r$-tuple of integers $\alpha = (\alpha_1,...,\alpha_r) \in \mathbb{Z}^r$ we have an endomorphism $[\alpha] \colon G^r \to G^r$ which on the $i$-th factor is just the $\alpha_i$-th power map. We also have the multiplication morphism $\Pi \colon G^r \to G$ which is defined by $\Pi(h_1,...,h_r) = h_1h_2...h_r$.

Now define for each $\alpha \in \mathbb{Z}^r$ a morphism $\phi_\alpha \colon C \to G$ by $\phi_\alpha = \Pi \circ [\alpha] \circ \iota$. By construction, we have $\phi_\alpha(c_e) = e \in G$ for every $\alpha \in \mathbb{Z}^r$. Hence the $\phi_\alpha \colon (C,c_e) \to (G,e)$ are morphisms of pointed varieties. Furthermore, we have that 
$$ \{ \phi_\alpha(c_g)~|~\alpha\in\mathbb{Z}^r\} = \{g_1^{\alpha_1}...g_r^{\alpha_r}~|~\alpha\in\mathbb{Z}^r\} $$
is a Zariski-dense subset of $G$. By Proposition \ref{prop:geomspeccriterion}, this implies that the $\phi_\alpha$ form a covering set through $e \in G$. As $G$ has a group structure, this covering set can now be translated to go through an arbitrary point.
\end{proof}

Over an uncountable base field, we may also use Lemma \ref{lemma:finitegeneratingset2} instead of Lemma \ref{lemma:finitegeneratingset} to construct covering sets which avoid any prespecified set of codimension at least two.

\begin{proposition}\label{prop:puncturedgroupschemegeomspec}
Let $k$ be an algebraically closed field of characteristic zero whose transcendence degree over $\mathbb{Q}$ is infinite.
Let $G$ be a connected algebraic group over $k$ and let $\Delta \subseteq G$ be a closed subset of codimension at least two.
Then $G \setminus \Delta$ is geometrically special.
\end{proposition}
\begin{proof}
By \cite[Lemma~2.16]{JRGeomSpec}, we may assume that $k$ is uncountable.
(Indeed, any covering set contains a countable subset which is still a covering set and similarly the set $S$ in Definition \ref{def:maindefinition} also always contains a countable Zariski-dense subset. Thus, we always only need to reference countably many transcendent variables when verifying geometric specialness.)

The set of $g \in G(k)$ such that $g \Delta$ does not contain the neutral element $e \in G$ is Zariski-dense in $G$. Thus, to show that $G \setminus \Delta$ is geometrically special, we may assume that $\Delta$ does not contain the neutral element and it suffices to construct a covering set through the origin of $G$ under this assumption.
Let $g_1,...,g_r \in G(k)$ be elements as in Lemma \ref{lemma:finitegeneratingset2}. Consider again the quasi-projective variety $G^r$. Let $\widetilde{\Delta} \subseteq G(k)^r$ be the set of all $(h_1,...,h_r) \in G(k)^r$ such that there is an $(\alpha_1,...,\alpha_r) \in \mathbb{Z}^r$ such that $h_1^{\alpha_1}...h_r^{\alpha_r} \in \Delta$. Note that $\widetilde{\Delta}$ is a countable union of closed subsets of codimension at least two in $G^r$. Let $C \to G^r$ be a very general curve connecting $(e,e,...,e)$ and $(g_1,...,g_r)$. Then $C$ lies entirely within $G^r \setminus \widetilde{\Delta}$. Now proceed exactly as in the proof of Proposition \ref{prop:groupschemegeomspec}.
\end{proof}

For the sake of completeness, let us prove that connected algebraic groups are Campana-special, since this appears to be missing from the literature.
To do so, we will make use of the following fibration lemma for Campana-special varieties.

\begin{lemma} \label{lemma:campana_special_fibration}
Let $f \colon X \to Y$ be a morphism of smooth varieties.
Assume that $f$ is surjective with connected fibers and has no nowhere reduced fibers.
Assume furthermore that both $Y$ and the (geometric) generic fiber of $f$ are Campana-special.
Then $X$ is Campana-special.
\end{lemma}
\begin{proof}
This follows from known properties of Campana-special varieties using Campana's notions of the \emph{core} of a smooth variety and the \emph{orbifold base} of a morphism of varieties (see \cite[Section~10]{CampanaJussieu} for more details, including the proof that the core exists).
Let us elaborate a bit.

Consider the core map $X \ratmap c(X)$ and let $\Delta$ denote its orbifold base divisor.
Then the pair $(c(X), \Delta)$ is of general type and in order to show that $X$ is Campana-special, it suffices to show that $c(X)$ is zero-dimensional.
Since the generic fiber of $f$ is Campana-special, the core $X \ratmap c(X)$ contracts the fibers of $f$.
Consequently, the map $X \ratmap c(X)$ factors over a dominant strictly rational map $Y \ratmap c(X)$.
Replacing $X$ and $Y$ by a blowup if necessary, we may assume that $X \to c(X)$ and $Y \to c(X)$ are morphisms.
We now claim that the orbifold base of the morphism $Y \to c(X)$ is also given by $\Delta$.
To see this, let $\eta \in c(X)$ be a point of codimension one and let $X_\eta$ and $Y_\eta$ denote the respective scheme-theoretic fibers over $\eta$.
Then $X_\eta$ and $Y_\eta$ are a finite set of codimension one points of $X$ and $Y$, respectively and moreover, $X_\eta$ is the preimage of $Y_\eta$ in $X$.
Since $X \to Y$ has no nowhere reduced fibers in codimension one, the minimum over the multiplicities of the points appearing in $X_\eta$ is the same as the minimum over the multiplicities of the points appearing in $Y_\eta$.
In other words, the coefficient of $\overline{ \{ \eta \} }$ in $\Delta$ is the same as the coefficient of $\overline{ \{ \eta \} }$ in the orbifold base of $Y \to c(X)$, which implies our claim.
As $Y$ is Campana-special, the orbifold base of any morphism from $Y$ to a positive-dimensional variety cannot be of general type. 
Hence, the variety $c(X)$ must be zero-dimensional.
Consequently, we see that $X$ is indeed Campana-special, as desired.
\end{proof}

\begin{lemma} \label{lemma:alg_groups_special}
Let $G$ be a connected algebraic group over $k$.
Then $G$ is Campana-special.
\end{lemma}
\begin{proof}
By Chevalley's structure theorem \cite[Theorem~1.1]{ConradChevalley}, $G$ admits a surjection $\pi \colon G \to A$ onto an abelian variety $A$ such that the fibers of $\pi$ are connected linear algebraic groups.
Since every smooth projective variety of Kodaira dimension zero is Campana-special \cite[Theorem~5.1]{CampanaFourier}, we see that $A$ is Campana-special.
(Alternatively, since every sheaf $\Omega^p_A$ is a trivial vector bundle, one can also deduce this directly from the definition.)
Thus, combining the above with Lemma~\ref{lemma:campana_special_fibration}, it suffices to show that linear algebraic groups are Campana-special.
So let $H$ be a connected linear algebraic group.
Then, the union of all Cartan subgroups of $H$ contains a dense open subset \cite[Theorem~6.4.5]{SpringerAlgGroups}.
In particular, given two general points $x, y \in H$, we can find a Cartan subgroup of $H$ containing $xy^{-1}$.
Thus, there is a coset of a Cartan subgroup of $H$ containing both $x$ and $y$.
Since Cartan subgroups are nilpotent, they are isomorphic to the product of a maximal torus of $H$ with a unipotent group \cite[Proposition~6.4.2 and Corollary~6.3.2]{SpringerAlgGroups}.
Over a perfect field, the underlying variety of a unipotent group is always isomorphic to affine space \cite[Corollary~14.3.10 and Corollary~14.2.7]{SpringerAlgGroups}.
Hence we see that the underlying variety of a Cartan subgroup is given by $(\mathbb{A}^1 \setminus \{0\})^m \times \mathbb{A}^n$ for some integers $m$ and $n$.
Since $\mathbb{A}^1 \setminus \{0\}$ and $\mathbb{A}^1$ are clearly Campana-special, an iterated application of Lemma \ref{lemma:campana_special_fibration} thus shows that Cartan subgroups are always Campana-special.
Consequently, we have shown that any two general points of $H$ can be joined by a Campana-special subvariety of $H$.
(In fact, the same argument shows that two arbitrary points $x,y \in H$ can be joined by a \emph{chain} of Campana-special subvarieties:
Simply choose a $z \in H$ such that both $xz^{-1}$ and $yz^{-1}$ are contained in some Cartan subgroup; not necessarily the same one.)
It now follows from \cite[Corollaire~10.8]{CampanaJussieu} that this implies that $H$ is Campana-special as well, finishing the proof.
\end{proof}

\section{Elliptic K3 surfaces and Enriques surfaces}

In this section we show that elliptic K3 surfaces are always geometrically special.
As a direct consequence, we also see that all Enriques surfaces are geometrically special.
Furthermore, we show that this geometric specialness continues to hold after removing finitely many points, as long as these points are chosen generically enough.
Our line of attack is to use that elliptic K3 surfaces always admit a non-torsion multisection of genus $\leq 1$ (note that in the absence of a zero section, the word ``non-torsion'' is somewhat ambiguous; see Definition \ref{def:nt_multisection} below for the precise definition).
This was originally established by Bogomolov and Tschinkel to prove potential density of rational points on elliptic K3 surfaces \cite{BogTschK3} -- however beware that one of the lemmas in loc.\ cit.\ is wrong for isotrivial fibrations, as is pointed out in \cite[Section~6]{Hassett2003}.
It is for this reason that we instead use a modified version of their claim due to Hassett and Kollár (see Lemma \ref{lemma:exists_non_torsion_multisection} below).
As our methods apply to all elliptic surfaces admitting such a multisection, the section will be written in that generality.

\subsection{Preliminaries on elliptic surfaces}

Since our approach fundamentally relies on the theory of elliptic fibrations on surfaces, let us briefly summarize the parts of the theory that we will use here.
The results below are all well-known, we refer to \cite{BTESMiranda} for proofs.

An \emph{elliptic fibration} is a proper morphism $\pi \colon X \to B$ from a smooth projective surface onto a smooth projective curve whose generic fiber $X_\eta$ is a smooth curve of genus one.
The elliptic fibration $\pi$ is said to be \emph{minimal} if no fiber of $\pi$ contains a $(-1)$-curve (i.e. a smooth curve of genus zero whose self-intersection is $-1$).
As $(-1)$-curves can always be contracted, every elliptic fibration can be blown down to a minimal elliptic fibration, which is unique.

An elliptic fibration $\pi \colon X \to B$ is said to be \emph{Jacobian} if $\pi$ admits a section $Z \colon B \to X$.
The choice of such a section yields a marked point of every fiber of $\pi$, turning the smooth fibers into elliptic curves.
In this way, the smooth fibers of $\pi$ become algebraic groups (note that this includes the generic fiber $X_\eta$).
In particular, this allows us to add any two sections $\sigma_1, \sigma_2 \colon B \to X$ by first adding them on the generic fiber and then spreading out the point of $X_\eta$ thus obtained to a section of $\pi$.
This induces a group structure on the set of sections of $\pi$ whose neutral element is $Z$; this group is called the \emph{Mordell-Weil group} and is denoted $\MW(X_\eta)$.
Moreover, we obtain for every integer $n \in \mathbb{Z}$ a rational map $[n] \colon X \ratmap X$ preserving $\pi$ which on the smooth fibers is just multiplication by $n$.
One can show that $[n]$ induces an endomorphism $[n] \colon X^o \to X^o$ where $X^o \subseteq X$ denotes the open subset on which $\pi$ is smooth.

\subsection{Torsion specializations of non-torsion multisections}

If the elliptic fibration $X \to B$ is Jacobian, i.e., if it admits a section which we can use as the neutral element to get a fibrewise group structure, our approach to proving geometric specialness can be summarized as follows:
Take any rational non-torsion multisection $C \to X$ and repeatedly add it to itself using the fibrewise group law.
As the multisection is non-torsion, this produces infinitely many pairwise distinct rational multisections, which therefore cover the surface.
While the graphs of these multisections are not dense in the product space (they are all contained in $X \times_B C \subseteq X \times C$), the rationality of the multisection means that we can reparametrize it in infinitely many different ways to solve this issue.
The other problem is that to obtain a covering set, we need to ensure that all of these multisections pass through the same point on the elliptic surface.
More precisely, to prove geometric specialness, we need to show that the set of points through which infinitely many of these multisections pass is Zariski-dense in the elliptic surface.
The goal of this subsection is prove that this is indeed the case.
More precisely, we will show the following proposition:

\begin{proposition}\label{prop:inftorsi_primep}
Let $\pi \colon X \to B$ be a Jacobian elliptic fibration with zero section $Z \colon B \to X$.
Assume that $\chi(X, \mathcal{O}_X) > 0$.
Let $\sigma \colon C \to X$ be a non-torsion multisection.
Then, for all but finitely many prime numbers $p$, there is a closed point $c \in C$ such that $\sigma(c)$ has order $p$ in its fiber.
\end{proposition}

Our proof relies on an intersection theoretic calculation, which also gives the estimate that the number of such points $c \in C$ grows quadratically in $p$.
We note that if the elliptic fibration $X \to B$ is not isotrivial, i.e. if there is no finite covering $B'$ such that the fiber product $X \times_B B'$ is birational to the product $E \times B'$ for some elliptic curve $E$, then we can substantially shorten the proof of Proposition~\ref{prop:inftorsi_primep} by appealing to Siegel's theorem over function fields; see Remark~\ref{rem:siegel} below.

\begin{lemma}\label{lemma:nointersect}
Let $\pi \colon X \to B$ be a minimal Jacobian elliptic fibration with zero section $Z \colon B \to X$.
Let $C_p \subseteq X$ be the closure of the curve containing all the points of precise order $p$ for some prime number $p$.
Then $C_p$ and $Z$ have no intersection points.
\end{lemma}
\begin{proof}
Let $K$ be the function field of $B$.
Then the generic fiber $X_\eta$ of $\pi$ is an elliptic curve over $K$.
Choose a finite algebraic field extension $K\subseteq L$ such that all $p$-torsion of $X_\eta$ is $L$-rational and let $C \to B$ be the morphism of smooth curves that corresponds to $K\subseteq L$.
Consider the pullback $\pi_Y \colon Y \to C$ of $\pi$ along $C\to B$.
By construction, $\pi_Y$ is a Jacobian elliptic fibration and the generic fiber $Y_\eta$ has $L$-rational $p$-torsion.
Let $X^o \subseteq X$ denote the open subset where $\pi$ is smooth.
Since $Z$ is a section, its image is contained in $X^o$.
As being smooth is stable under pullbacks, $\pi_Y$ is smooth on the open subset $X^o \times_B C \subseteq Y$ and consequently, $X^o \times_B C$ is a smooth variety.
Choose a resolution of singularities $Y' \to Y$ that is an isomorphism over $X^o \times_B C$ and blow down any vertical $(-1)$-curves on $Y'$ to obtain a smooth surface $Y''$ equipped with a minimal elliptic fibration $\pi_Y^{min} \colon Y'' \to C$.
The surface $Y''$ is birational to $Y$ and thus, the generic fiber of $\pi_Y^{min}$ still has $L$-rational $p$-torsion.
As torsion sections on a minimal elliptic surface do not intersect (see \cite[Proposition~VII.3.2]{BTESMiranda}), the zero section intersects no nontrivial $p$-torsion section on $Y''$ and consequently also not on $Y'$, since a blowup cannot introduce new intersection points.
As the zero section does not leave $X^o \times_B C$, there still are no intersection points after the blowdown $Y' \to Y$.
Hence there are also no intersection points on $X$. 
\end{proof}

\begin{remark}
The above argument also shows that two different torsion multisections have no intersection points in $X^o$.
If we however do not assume one of the multisections to be an actual section, there can be intersection points outside of $X^o$.
For an explicit example consider the Weierstrass equation $y^2 = x(x-\lambda)(x-2\lambda)$ with parameter $\lambda\in\mathbb{A}^1$.
The equation degenerates to a cuspidal curve when $\lambda=0$.
Every non-zero torsion point of $y^2 = x(x-\lambda)(x-2\lambda)$ specializes to the cusp when $\lambda = 0$, so all non-trivial torsion multisections intersect there.
\end{remark}

\begin{lemma}\label{lemma:numeq}
Let $\pi \colon X \to B$, $Z$ and $C_p$ be as in the previous lemma. Then there is a numerical equivalence
$$ C_p \sim (p^2-1)Z + (p^2-1)\chi(\mathcal{O}_X) F + R $$
where $F$ is the class of a closed fiber of $\pi$ and $R$ is a linear combination of fiber components that do not meet $Z$.
\end{lemma}
\begin{proof}
Consider the canonical group homomorphism $\Num(X) \to \MW(X_\eta)$ given by intersecting a divisor with $X_\eta$.
The divisor $C_p$ is by definition all of the nontrivial $p$-torsion of $X_\eta$ and is consequently in the kernel of that group homomorphism.
This kernel is generated by the class of $Z$, $F$, and the classes of the irreducible components of the singular fibers \cite[Theorem~VII.2.1]{BTESMiranda}.
Thus, we have a numerical equivalence
$$ C_p \sim aZ + bF + R $$
for some $a,b \in \mathbb{Z}$ and some $R$ as in the statement.
It remains to determine the values of $a$ and $b$. \\
For this, we intersect both sides of the equivalence with $Z$ and $F$ and compute the resulting intersection numbers.
We have $Z.F = 1$ as $Z$ is a section and $Z.R = 0$ by definition of $R$.
Furthermore, we have $F.F = 0$ and $F.R = 0$.
As the elliptic fibration was assumed to be Jacobian, there are no multiple fibers.
Thus, the canonical class of $X$ is numerically equivalent to $(2g_B - 2 + \chi(\mathcal{O}_X))F$ \cite[Corollary~V.12.3]{CCS}.
By adjunction, it follows that $Z.Z = -\chi(\mathcal{O}_X)$.
Furthermore, $C_p.F = p^2-1$ as an elliptic curve over an algebraically closed field has $p^2-1$ nontrivial $p$-torsion points and $C_p.Z = 0$ by Lemma \ref{lemma:nointersect}.
Combining all of these computations we find $a=p^2-1$ and $b = (p^2-1)\chi(\mathcal{O}_X)$, as desired.
\end{proof}

\begin{proof}[Proof of Proposition~\ref{prop:inftorsi_primep}]
We reuse the notation $C_p$ from the previous lemmas.
We may assume that $\sigma \neq Z$ and that $\sigma$ and $C_p$ have no common components, as otherwise the statement is trivially true.
Let us first assume that $\sigma$ is a section of $\pi$.
In that case we have $B=C$ and $\sigma$ does not hit any of the singular points of the singular fibers.
Thus, if we replace $\sigma$ by $m\sigma$ for a sufficiently composite number $m \in \mathbb{N}$, it will only go through the identity component of every singular fiber.
Then, using Lemma \ref{lemma:numeq}, we have:
$$ C_p.\sigma = (p^2-1)Z.\sigma + (p^2-1)\chi(\mathcal{O}_X) $$
As we assumed $\chi(\mathcal{O}_X) > 0$, we see $C_p.\sigma > 0$.
This means that there is at least one $b \in B$ such that $\sigma(b)$ is of exact order $p$ in its fiber and this holds true for every prime number $p$.
A point in $X^{sm}$ cannot be simultaneously of precise order $p_1$ and $p_2$ if $p_1 \neq p_2$, so we are done. \bigbreak \par
In general, $\sigma$ will not be a section and in that case we have a finite ramified covering $\pi \circ \sigma \colon C \to B$.
Consider the following cartesian square:
\begin{equation*} \begin{tikzcd}
X \times_B C \ar[r] \ar[d] & X \ar[d, "\pi"] \\ C \ar[r, "\pi\circ\sigma"] & B
\end{tikzcd} \end{equation*}
The fiber product $X \times_B C$ will in general not be smooth, so we take a resolution of singularities $X' \to X \times_B C$.
The induced morphism $X' \to C$ is still an elliptic fibration.
If $X'$ has any vertical \mbox{$(-1)$-curves} we blow them down to get a minimal elliptic fibration $X'' \to C$.
This elliptic fibration now has two different sections, which are given by lifting the two sections $(Z\circ\pi\circ\sigma, \id_C) \colon C \to X \times_B C$ and $(\sigma, \id_C) \colon C \to X \times_B C$ to $X''$.
By the previous case we considered, the section $(\sigma, \id_C)$ (regarded as a section of $X'' \to C$) assumes torsion values in infinitely many fibers.
As the blowups only change finitely many fibers, the same is true on $X \times_B C$.
Hence the original multisection $\sigma \colon C \to X$ must have assumed infinitely many torsion values as well.
\end{proof}

\begin{remark}\label{rem:siegel}
If the fibration $X \to B$ is not isotrivial, then Proposition~\ref{prop:inftorsi_primep} admits a quicker, but less effective, proof using Siegel's theorem over function fields \cite[Theorem~III.12.1]{AdvEllCurves}.
(It states that given a non-constant elliptic curve $E$ over the coordinate ring $R$ of a smooth affine curve and an affine open $U \subseteq E$, the set $U(R)$ is finite.)

Let $Y$ be a desingularization of the fiber product $X \times_B C$.
Then the projection $\widetilde{\pi} \colon Y \to C$ is an elliptic fibration and both $Z$ and $\sigma$ induce sections $\widetilde{Z}, \widetilde{\sigma} \colon C \to Y$.
Treating $\widetilde{Z}$ as the zero section, the assumption that $\sigma$ is non-torsion implies that $\widetilde{\sigma}$ is a non-torsion section.
Replacing $C$ by a dense open if necessary, we may assume that the map $\widetilde{\pi}$ is smooth and that moreover, for every $c \in C$, we have $\widetilde{\sigma}(c) \neq \widetilde{Z}(c)$.
Consider the sections $p \widetilde{\sigma}$ as $p$ runs over all prime numbers.
Since $Y \setminus \widetilde{Z} \to C$ is an affine open of a relative elliptic curve, it has only finitely many sections by Siegel's theorem over function fields.
Thus, it follows that for almost all prime numbers $p$, the section $p \widetilde{\sigma}$ intersects the zero section $\widetilde{Z}$.
In particular, for almost all prime numbers $p$, there is a point $c \in C$ such that $p \widetilde{\sigma}(c)=0$ in its fiber.
Since $\widetilde{\sigma}(c)$ is not zero itself, it follows that $\widetilde{\sigma}(c)$ is torsion of exact order $p$. 
\end{remark}

\subsection{Geometric specialness of certain elliptic surfaces}

\begin{proposition}\label{prop:jacobianellfibgeomspec}
Let $\pi \colon X \to B$ be a Jacobian elliptic fibration with zero section $Z \colon B \to X$.
Assume that $\chi(X, \mathcal{O}_X) > 0$ and that $X$ admits a non-torsion multisection $\sigma \colon C \to X$ with $C$ smooth of genus $g_C \leq 1$.
Then $X$ is geometrically special.
\end{proposition}
\begin{proof}
Since $\chi(X,\mathcal{O}_X)>0$, by Proposition~\ref{prop:inftorsi_primep}, we see that for almost every prime number $p$, the multisection $\sigma$ assumes a torsion value of order $p$ in some fiber. Consequently, the following subset of $X$ is Zariski-dense:
$$ \{ m\sigma(c)~|~m\in\mathbb{N},~c\in C,~\sigma(c)~\text{is torsion in its fiber} \} $$
Let $s = m\sigma(c)$ be any point in this set and let $n$ be the order of $\sigma(c)$ in its fiber. We now construct a covering set for $X$ through the point $s$.

For this, we need an $\mathbb{N}$-indexed sequence of pairwise distinct endomorphisms $\phi_i \colon C \to C$ fixing the point $c$. If $g_C=0$, then we take any infinite sequence of automorphisms fixing $c$. If $g_C=1$, we take $\phi_i$ to be multiplication by $i$ with respect to the origin $c$. These morphisms have the property that the set $\{ \phi_i(d)~|~i \in \mathbb{N} \}$ is infinite for infinitely many $d \in C$: This is clear for $g_C=0$, and for $g_C=1$ it follows from the existence of infinitely many non-torsion points on $C$.

A covering set is then given by the following morphisms, with $i \in \mathbb{N}, k \in \mathbb{N}$:
$$ f_{ik} \colon C \to X \quad f_{ik}(x) = ([kn+m] \circ \sigma)(\phi_i(x)) $$
Note that while the multiplication-by-$(kn+m)$ map is only defined over a dense open of $X$, the normality of $C$ and the properness of $X$ ensure that the $f_{ik}$ do in fact define morphisms.
By construction, $f_{ik}(c)=s$ for all $i,k$. For the density of the graphs consider a point $d\in C$ such that $\{ \phi_i(d)~|~i \in \mathbb{N} \}$ is infinite. Then, as $\sigma$ is non-torsion, the order of $\sigma(\phi_i(d))$ in its fiber is unbounded when varying $i$. Consequently the values of $f_{ik}(d)$ are dense in $X$ when varying $i,k$ (but fixing $d$). Now recall that infinitely many such $d$ exist, so the union of the graphs of the $f_{ik}$ is dense in $C \times X$. (Alternatively, a single such $d$ is enough, by Proposition \ref{prop:geomspeccriterion}.)
\end{proof}

The previous proposition only applies to Jacobian elliptic fibrations, as we crucially used the existence of a multiplication-by-$n$ map to move our fixed non-torsion multisection around.
However, as our intended application are K3 surfaces and the elliptic fibration on a K3 surface need not be Jacobian, we need to extend our methods to also cover the non-Jacobian case.
To do so, we first need to clarify what a non-torsion multisection is in the absence of a zero section.
We will use the following notion, which was already used by Bogomolov--Tschinkel in \cite[Definition~3.7]{BogTschK3} under the name ``nt-multisection''.

\begin{definition}\label{def:nt_multisection}
Let $\pi \colon X \to B$ be an elliptic fibration.
A multisection $\sigma \colon C \to X$ is said to have \emph{torsion differences of order $m$} if for every $c_1, c_2 \in C$ such that $\pi(\sigma(c_1)) = \pi(\sigma(c_2))$, the divisor $\sigma(c_1)-\sigma(c_2)$ is torsion of order $m$ in $\Pic^0(X_{\pi(\sigma(c_1))})$.
If there is no $m$ such that $\sigma$ has torsion differences of order $m$, we say that $\sigma$ has \emph{non-torsion differences}.
\end{definition}

Let us recall some general facts about non-Jacobian elliptic fibrations.
Let $\pi \colon X \to B$ be a minimal elliptic fibration, which we will assume to have no multiple fibers (this is purely for convenience; but this always holds in the case of elliptic K3 surfaces).
We can then find an associated Jacobian elliptic fibration $\pi_J \colon J \to B$ over the same base, where the generic fiber $J_\eta \subseteq J$ consists of divisors of degree $0$ on the generic fiber $X_\eta \subseteq X$.
Consequently, $J_\eta$ contains the empty divisor as a distinguished element and this element corresponds to the zero section $Z \colon B \to J$.
We call $J$ the \emph{Jacobian fibration} of $X$ or the \emph{Jacobian} for short.
Note that $J_\eta$ is simply the Jacobian curve of the smooth genus one curve $X_\eta$.
Also note that if the fibration $X \to B$ is already Jacobian, it is isomorphic to its Jacobian fibration.
When $b \in B$ is a closed point, the fibers $X_b$ and $J_b$ are isomorphic (this is only true in the absence of multiple fibers).
In particular it follows that $X \to B$ and $J \to B$ have the same number of singular fibers and that $\chi(\mathcal{O}_X) = \chi(\mathcal{O}_J)$.
We refer to \cite[Chapter V]{CossecDolgachev} for details.

By construction, the generic fiber $X_\eta$ is a torsor under the elliptic curve $J_\eta$.
Thus there is a morphism $J_\eta \times_{\kappa(\eta)} X_\eta \to X_\eta$, which defines a rational map $J \times_B X \ratmap X$.
It can be shown that this map is defined on every smooth fiber of $X$.
When $b \in B$ is a point over which $\pi$ is smooth, this rational map restricts to a transitive group action of $J_b$ on $X_b$.
In this way, the $B$-scheme $X$ is a torsor under the group $B$-scheme $J$.

Suppose that there is an integral curve $D \subseteq X$ such that the induced morphism $D \to B$ is of degree $d$.
Then $D_\eta$ is a divisor of degree $d$ on the generic fiber $X_\eta$.
Given such a curve $D$, there exists a rational map $\eta_D \colon X \ratmap J$ over $B$ which is finite of degree $d^2$ when restricted to the smooth fibers.
The map $\eta_D$ is easy to describe on the Picard groups:
It maps divisors of degree $1$ to divisors of degree $0$ by first multiplying them by $d$ and then subtracting the divisor $D$.
This gives a well-defined morphism over the smooth fibers of $X$, i.e.\ a rational map $X \ratmap J$, and it preserves the fibers by construction.
The map $\eta_D$ has the following useful property:

\begin{lemma}
Let $\sigma \colon C \to X$ be a multisection with non-torsion differences.
Then $\eta_D \circ \sigma \colon C \to J$ is a multisection of the same degree, has non-torsion differences, and the induced rational map $\sigma(C) \ratmap \eta_D(\sigma(C))$ is birational.
\end{lemma}
\begin{proof}
See \cite[Lemma~3.9]{BogTschK3}.
\end{proof}

Thus, the map $\eta_D$ allows us to transport multisections with non-torsion differences from an elliptic fibration to its Jacobian.
We now use this property and the $J$-torsor structure of $X$ to prove an analogous result to Proposition \ref{prop:jacobianellfibgeomspec} for the case of a non-Jacobian elliptic fibration.

\begin{proposition}\label{prop:ellfibgeomspec}
Let $\pi \colon X \to B$ be a non-Jacobian elliptic fibration without multiple fibers.
Assume that $\chi(X, \mathcal{O}_X) > 0$ and that $X$ admits a multisection $\sigma \colon C \to X$ with non-torsion differences, where $C$ is smooth of genus $g_C \leq 1$.
Then $X$ is geometrically special.
\end{proposition}
\begin{proof}
Let $\pi_J \colon J \to B$ be the associated Jacobian elliptic fibration, and let $\eta_C \colon X \ratmap J$ be the rational map described above.
Then $\eta_C \circ \sigma \colon C \to J$ is a multisection with non-torsion differences and the induced maps $\pi_J \circ \eta_C \circ \sigma$ and $\pi \circ \sigma \colon C \to B$ agree.
Consequently, we can use the $J$-torsor structure of $X$ to move $\sigma$ around by the multiples of $\eta_C \circ \sigma$.
Since we have $\chi(X, \mathcal{O}_X) > 0$, we also have $\chi(J, \mathcal{O}_J)>0$.
Thus, Proposition \ref{prop:inftorsi_primep} applies and the multisection $\eta_C \circ \sigma$ assumes infinitely many torsion values.
As $\eta_C \circ \sigma$ is non-torsion, the torsion orders are unbounded.
Hence, the set
$$ \{ (m~\eta_C(\sigma(c)))\cdot \sigma(c)~|~m\in\mathbb{N},~c \in C,~\eta_C(\sigma(c))~\text{is torsion in its fiber} \} $$
is dense in $X$, where the $\cdot$ denotes the $J$-action on $X$.
Let $s = (m~\eta_C(\sigma(c)))\cdot \sigma(c)$ be any point in this set and let $n$ be the order of $\eta_C(\sigma(c))$ in its fiber.
We now construct a covering set for $X$ through the point $s$.

Choose an $\mathbb{N}$-indexed sequence of pairwise distinct endomorphisms $\phi_i \colon C\to C$ fixing $c$.
There are infinitely many points $d \in C$ such that $\phi_i(d)$ takes on infinitely many values as $i$ varies.

A covering set is then given by the following morphisms, with $i,k\in\mathbb{N}$:
$$ f_{ik} \colon C \to X \quad f_{ik}(x) = ((nm+k)~\eta_C(\sigma(\phi_i(x)))) \cdot \sigma(\phi_i(x)) $$
Note that the $J$-torsor structure of $X$ and the multiplication-by-$(nm+k)$ map are only defined over a dense open subset of $B$, so the $f_{ik}$ are a priori only rational maps.
As $C$ is a smooth curve and $X$ is proper, they extend to morphisms.
The density of the graphs of the $f_{ik}$ follows just as in Proposition \ref{prop:jacobianellfibgeomspec}, recalling that the action of $J$ on $X$ is transitive on smooth fibers. 
\end{proof}

It appears to be an open problem whether every elliptic fibration $X \to B$ with base curve of genus $g_B \leq 1$ admits a multisection $C \to X$ of genus $g_C \leq 1$ with non-torsion differences.
In the particular case of elliptic K3 surfaces however, this is always the case.

\begin{lemma}\label{lemma:exists_non_torsion_multisection}
Let $X$ be an elliptic K3 surface and let $\pi \colon X \to \mathbb{P}^1$ be an elliptic fibration. Then there exist an elliptic curve $E$ and a multisection $E \to X$ with non-torsion differences.
\end{lemma}
\begin{proof}
See \cite[Corollary~10.9]{Hassett2003}.
\end{proof}

\begin{remark}
For the sake of completeness, we note the following result due to Baltes \cite[Theorem~4.2]{BaltesK3}:
If $X$ is an elliptic K3 surface with elliptic fibration $\pi \colon X \to \mathbb{P}^1$, then there is a (potentially different) elliptic fibration $\pi' \colon X \to \mathbb{P}^1$ which admits a \emph{rational} multisection with non-torsion differences. 
Moreover, by \cite[Proposition~9.6]{Hassett2003}, one can choose $\pi = \pi'$ except possibly when $X$ is a Kummer surface or when $X$ has Picard rank $20$.
\end{remark}

\begin{corollary}\label{cor:ellk3geomspec}
Let $X$ be an elliptic K3 surface. Then $X$ is geometrically special.
\end{corollary}
\begin{proof}
The elliptic fibration of a K3 surface automatically has no multiple fibers \cite[Proposition~11.1.6]{HuybrechtsK3}.
By Lemma \ref{lemma:exists_non_torsion_multisection}, the K3 surface $X$ admits an elliptic fibration with a multisection of genus $\leq 1$ with non-torsion differences.
Thus, one of the two Propositions \ref{prop:ellfibgeomspec} or \ref{prop:jacobianellfibgeomspec} applies and we conclude.
\end{proof}

\begin{corollary}\label{cor:enriquesgeomspec}
Let $X$ be an Enriques surface. Then $X$ is geometrically special.
\end{corollary}
\begin{proof}
Every Enriques surface in characteristic zero admits a 2-to-1 étale cover $Y \to X$ where $Y$ is a K3 surface.
As every Enriques surface admits an elliptic fibration $X \to \mathbb{P}^1$, we obtain a surjective morphism $Y \to \mathbb{P}^1$.
Every morphism $Y \to \mathbb{P}^1$ is automatically an elliptic fibration \cite[Remark~11.1.5.(iv)]{HuybrechtsK3}.
Hence $Y$ is geometrically special by Corollary \ref{cor:ellk3geomspec} and thus, so is $X$. 
\end{proof}

\subsection{Puncturing elliptic surfaces}

In this section, we show how to modify the arguments of the previous section to deal with the case of elliptic K3 surfaces deprived of finitely many closed points.
Our results can be summarized by saying that as long as the points are chosen generically enough, the arguments still go through.
As a consequence, we see that the complement of a generic set of finitely many points in an elliptic K3 is geometrically special.

We will need the following elementary lemma about congruences of integers.

\begin{lemma}\label{lemma:elementarycongruences}
Let $\delta \geq 0$, $n \geq 1$, $a_1,...,a_\delta \in \mathbb{Z}$ and $b_1,...,b_\delta > \delta$ be integers. Assume that $n$ is coprime to all the $b_i$. Then, for any integer $m$, there are infinitely many integers $k \in \mathbb{Z}$ such that $k \equiv m \mod n$ and $k \not \equiv a_i \mod b_i$ for $i=1,...,\delta$.
\end{lemma}
\begin{proof}
We write $B := b_1...b_\delta$. It suffices to show that there is a single integer $k$ satisfying these congruences, as if $k$ is a solution, then so is $k + n B$. Now observe that of the $B$ integers between $0$ and $B-1$, for each $i=1,...,\delta$, exactly $B/b_i$ satisfy $k \equiv a_i \mod b_i$. As each of the $b_i$ is strictly greater than $\delta$, we have $B/b_1 + ... + B/b_\delta < B$. Thus, there is an integer $k$ satisfying none of the congruences $k \equiv a_i \mod b_i$, i.e. $k \not \equiv a_i \mod b_i$ for every $i=1,...,\delta$. As $n$ was assumed to be coprime to all $b_i$, by the Chinese Remainder Theorem, the integer $k$ can additionally be chosen to satisfy $k \equiv m \mod n$.
\end{proof}

Our techniques only allow us to remove finite sets of points which are in general position.
Indeed, we do not have much control over the torsion specializations of the multisections that we use to verify geometric specialness.
In particular, we cannot exclude the case that there is a Jacobian elliptic K3 surface for which \emph{all} multisections of genus $\leq 1$ pass through the same point.
For such a surface, removing that one point from the surface would fundamentally break the proof strategy of the previous subsection.
Of course, similar issues arise if there is a finite subset which intersects every multisection of genus $\leq 1$. 
Let us write down the precise technical condition that the points need to satisfy.

\begin{definition}
Let $\pi \colon X \to B$ be a Jacobian elliptic fibration with zero section $Z \colon B \to X$ and let $\sigma \colon C \to X$ be a non-torsion multisection.
Let $\delta \geq 0$ be an integer.
We say that a subset $\Delta \subseteq X(k)$ of cardinality $\delta$ is in \emph{sufficiently general position with respect to $\sigma$} if $\pi$ is smooth at every point of $\Delta$, no point of $\Delta$ is contained in a fiber for which $\sigma$ intersects the non-smooth locus of $\pi$ and moreover, $\Delta$ does not contain a point $x$ satisfying the following:
There is an integer $a \leq \delta \deg(\pi \circ \sigma)$ such that $x$ is an $a$-torsion point in its fiber and such that $\sigma$ intersects the fiber of $x$ in an $a$-torsion point.
\end{definition}

We note that for a fixed choice of $X$, $\pi$, $Z$ and $\sigma$, the set of those $\Delta \subseteq X(k)$ which are in sufficiently general position with respect to $\sigma$ form a nonempty, hence dense, open subset of the product space $X^{\delta}$.
This justifies the terminology ``sufficiently general''.
We now show that the covering sets constructed in Proposition~\ref{prop:ellfibgeomspec} have an infinite sub-covering set disjoint of any prespecified $\Delta \subseteq X$ in sufficiently general position.

\begin{lemma}\label{lemma:multisectionsavoidgenericpoints2}
Let $\pi \colon X \to B$ be an elliptic fibration with associated Jacobian fibration $\pi_J \colon J \to B$.
Assume $\chi(X, \mathcal{O}_X) > 0$.
Let $\sigma \colon C \to X$ be a multisection with non-torsion differences.
Let $\delta \geq 0$ be an integer and let $\Delta \subseteq X$ be a finite set of closed points of cardinality $\delta$ such that $\eta_C$ is defined at all points of $\Delta$ and such that $\eta_C(\Delta) \subseteq J$ is in sufficiently general position with respect to $\eta_C \circ \sigma \colon C \to J$.
Let $c_0 \in C$ be a point such that $\eta_C(\sigma(c_0))$ is torsion of order $n$ in its fiber and assume that $n$ is coprime to the order of every torsion point of the form $\eta_C(\sigma(c))$ with $\pi(c) \in \pi(\Delta)$.
Then, for any integer $m$, there are infinitely many integers $k \in \mathbb{Z}$ such that $(([k] \circ \eta_C \circ \sigma) \cdot \sigma)(c_0) = (([m] \circ \eta_C \circ \sigma) \cdot \sigma)(c_0)$ and such that the image of $([k] \circ \eta_C \circ \sigma) \cdot \sigma \colon C \to X$ is disjoint from $\Delta$.
(Here, $\cdot$ denotes the action of $J$ on $X$.)
\end{lemma}
\begin{proof}
Let $k$ be an integer, let $c \in C$ and let $x \in \Delta$ be a point.
Then $(([k] \circ \eta_C \circ \sigma) \cdot \sigma)(c) = x$ if and only if $x-\sigma(c) = ([k] \circ \eta_C \circ \sigma)(c)$.
Consequently, fixing $c$ and $x$, this happens only for multiple values $k_1, k_2$ if $(\eta_C \circ \sigma)(c)$ is a torsion point of order $a$, where $a$ is an integer dividing $k_2 - k_1$.
Moreover, in that situation, $x-\sigma(c)$ must also be a torsion point of order dividing $a$.
(Note that the ``sufficiently general position'' assumption implies that $\pi$ is smooth at both $x$ and $\sigma(c)$, so that this difference is well-defined.)
Now consider the following set:
\[ \{ (c,x) \in C \times \Delta~|~\pi(\sigma(c)) = \pi(x)~\text{and}~(\eta_C \circ \sigma)(c)~\text{is torsion in its fiber} \} \]
Note that this is evidently a finite set of cardinality $\leq \delta \deg(\pi \circ \sigma)$.
For each $(c,x)$ in this set, define the positive integer $b_{(c,x)}$ as the order of $\eta_C(\sigma(c))$ in its fiber.
Now suppose that $b_{(c,x)} \leq \delta \deg(\pi \circ \sigma)$.
Then the ``sufficiently general position'' assumption implies that $\eta_C(x)$ is not $b_{(c,x)}$-torsion in its fiber.
As we have $\eta_C(x)-\eta_C(\sigma(c)) = \deg(\pi \circ \sigma)(x-\sigma(c))$ and the left hand side is not $b_{(c,x)}$-torsion, we see that $x-\sigma(c)$ is not $b_{(c,x)}$-torsion either.
Consequently, $b_{(c,x)} \leq \delta \deg(\pi \circ \sigma)$ can only happen if $x$ lies on at most one of the multisections $([k] \circ \eta_C \circ \sigma) \cdot \sigma$.
Consider the remaining $b_{(c,x)}$, which now all satisfy $b_{(c,x)} > \delta \deg(\pi \circ \sigma)$.
For each of them, let $a_{(c,x)}$ be an integer such that $(([a_{(c,x)}] \circ \eta_C \circ \sigma) \cdot \sigma)(c) = x$ if it exists and discard $b_{(c,x)}$ if it does not.
By Lemma \ref{lemma:elementarycongruences}, there are infinitely many integers $k$ for which $k \not \equiv a_{(c,x)} \mod b_{(c,x)}$ for every remaining $(c,x)$ as above and also $k \equiv m \mod n$ (with $m$ and $n$ as in the lemma statement).
The condition $k \equiv m \mod n$ ensures that $(([k] \circ \eta_C \circ \sigma) \cdot \sigma)(c_0) = (([m] \circ \eta_C \circ \sigma) \cdot \sigma)(c_0)$ and moreover, for all but finitely many of these values of $k$, we have that $([k] \circ \eta_C \circ \sigma) \cdot \sigma \colon C \to X$ is disjoint from $\Delta$ by construction.
\end{proof}

\begin{proposition}\label{prop:ellfibgeomspecpunctured}
Let $\pi \colon X \to B$ be an elliptic fibration without multiple fibers.
Assume that $\chi(X, \mathcal{O}_X) > 0$ and that $X$ admits a multisection $\sigma \colon C \to X$ with non-torsion differences where $C$ is smooth of genus $g_C \leq 1$.
Let $\delta \geq 0$ be an integer and let $\Delta \subseteq X$ be a set of closed points of cardinality $\delta$ such that $\eta_C$ is defined at all points of $\Delta$ and such that $\eta_C(\Delta) \subseteq J$ is in sufficiently general position with respect to $\eta_C \circ \sigma$.
Then $X \setminus \Delta$ is geometrically special.
\end{proposition}
\begin{proof}
Let $\pi_J \colon J \to B$ be the associated Jacobian elliptic fibration, and let $\eta_C \colon X \to J$ be as above.
Let $n_0$ be the least common multiple of the torsion orders of the points in $\eta_C(\sigma(c))$, where $c$ ranges over the points of $C$ satisfying $\pi(c) \in \pi(\Delta)$.
Since $\chi(X, \mathcal{O}_X)>0$, Proposition~\ref{prop:inftorsi_primep} implies that for almost every prime number $p$, the multisection $\eta_C \circ \sigma$ of $\pi_J$ assumes a torsion value of exact order $p$ in some fiber.
Consequently, the following subset of $X$ is Zariski-dense:
\[ \{ (([m] \circ \eta_C \circ \sigma) \cdot \sigma)(c)~|~m \in \mathbb{N},~c\in C,~\eta_C(\sigma(c))~\text{is torsion of order coprime to}~n_0~\text{in its fiber} \} \]
Here, as before, the $\cdot$ denotes the action of $J$ on $X$ over $B$.
Let $s = (([m] \circ \eta_C \circ \sigma) \cdot \sigma)(c_0)$ be any point in this set and let $n$ be the order of $\eta_C(\sigma(c_0))$ in its fiber. 
We now construct a covering set for $X \setminus \Delta$ through the point $s$.

By Lemma~\ref{lemma:multisectionsavoidgenericpoints2}, there are infinitely many integers $k$ such that we have $(([k] \circ \eta_C \circ \sigma) \cdot \sigma)(c_0) = s$ and such that the image of the morphism $([k] \circ \eta_C \circ \sigma) \cdot \sigma \colon C \to X$ is disjoint from $\Delta$.
Let $K \subseteq \mathbb{Z}$ denote the set of such integers.
As in the proof of Proposition~\ref{prop:jacobianellfibgeomspec}, pick an $\mathbb{N}$-indexed sequence of pairwise distinct endomorphisms $\phi_i \colon C \to C$ fixing the point $c_0$.

A covering set for $X \setminus \Delta$ through $s$ is then given by the following morphisms, with $i \in \mathbb{N}, k \in K$:
\[ f_{ik} \colon C \to X \setminus \Delta \quad f_{ik}(c) = (([k] \circ \eta_C \circ \sigma) \cdot \sigma)(\phi_i(c)) \]
While a priori, the $f_{ik}$ define only rational maps $C \ratmap X$, the normality of $C$ and the properness of $X$ ensure that the $f_{ik}$ extend to morphisms $C \to X$, whose image is disjoint from $\Delta$ by construction.
That $f_{ik}(c_0) = s$ follows from the choice of $K \subseteq \mathbb{Z}$ and the $\phi_i$.
The density of the graphs of the morphisms in $C \times (X \setminus \Delta)$ can be checked on $C \times X$, where the arguments are just like in the proof of Proposition~\ref{prop:jacobianellfibgeomspec}, recalling that the action of $J$ on $X$ is transitive on the fibers.
\end{proof}

To state the following two corollaries, recall that we say that a property is satisfied for the \emph{general subset} $\Delta \subseteq X(k)$ of cardinality $\delta$ if there is a dense open subset of the product space $X^{\delta} = X \times \cdots \times X$ such that the property is satisfied for every element of this open subset. 

\begin{corollary}\label{cor:punctured_k3_geomspec}
Let $X$ be an elliptic K3 surface over $k$, let $\delta \geq 0$ be an integer and let $\Delta \subseteq X(k)$ be a general subset of cardinality $\delta$. Then $X \setminus \Delta$ is geometrically special.
\end{corollary}
\begin{proof}
This follows directly from Proposition~\ref{prop:ellfibgeomspecpunctured} and the arguments used to prove Corollary~\ref{cor:ellk3geomspec}.
\end{proof}

\begin{corollary}\label{cor:punctured_enriques_geomspec}
Let $X$ be an Enriques surface over $k$, let $\delta \geq 0$ be an integer and let $\Delta \subseteq X(k)$ be a general subset of cardinality $\delta$. Then $X \setminus \Delta$ is geometrically special.
\end{corollary}
\begin{proof}
This follows directly from Corollary \ref{cor:punctured_k3_geomspec} and the arguments used to prove Corollary~\ref{cor:enriquesgeomspec}.
\end{proof}

\bibliography{geomspec}{}
\bibliographystyle{alpha}
\end{document}